\documentclass[12pt]{article}
\usepackage{enumerate,amsmath,amssymb,amsthm}
\usepackage[left=2.2cm,bottom=2.7cm,top=2.7cm,right=2.2cm]{geometry}
\newcounter{alphthm}
\newtheorem{thm}{Theorem}

\newtheorem{cor}{Corollary}
\newtheorem{lem}{Lemma}

\newcommand{\be}{\begin{equation}}
\newcommand{\ee}{\end{equation}}
\newcommand{\ben}{\begin{enumerate}}
\newcommand{\een}{\end{enumerate}}
\newcommand{\pa}{{\partial}}

\newcommand{\pxi}{{\pa \over \pa x^i}}

\def\beq{\begin{equation}}
\def\eeq{\end{equation}}

\title{On m-th root Finsler metrics with special curvature properties}

\author{A. Tayebi and B. Najafi}
\begin{document}
\maketitle

\begin{abstract}
In this Note, we prove that every m-th root Finsler metric with isotropic Landsberg
curvature reduces to a Landsberg metric. Then, we show that every m-th root metric with
almost vanishing H-curvature has vanishing H-curvature.\footnote{ To cite this article: A. Tayebi and B. Najafi, On m-th Root Finsler Metrics with Special Curvature Properties, C. R. Acad. Sci. Paris, Ser. I 349 (2011) 691–-693.}
\end{abstract}

\section{Introduction}\label{Introduction}
Let $M$ be an $n$-dimensional $ C^\infty$ manifold. Denote  by $TM=\cup _{x \in M} T_x M $ the tangent space of $M$. Let $TM_{0} = TM \setminus \{ 0 \}$. Let $F=\sqrt[m]{A}$ be a Finsler metric on $M$, where $A$ is given by
\begin{equation}
A:=a_{i_{1}\dots i_{m}}(x)y^{i_{1}}y^{i_{2}}\dots y^{i_{m}}
\label{1.1}
\end{equation}
with $a_{i_{1}\dots i_{m}}$ symmetric in all its indices \cite{shLi1}. Then $F$ is called an $m$-th root Finsler metric. Let  $F$ be an $m$-th root Finsler metric on an open subset $U\subset \mathbb{R}^n$. Put
\[
A_i=\frac{\pa A}{\pa y^i}, \ \ \  A_{ij}=\frac{\pa^2 A}{\pa y^i\pa y^j},\ \ \ A_{x^k}=\frac{\pa A}{\pa x^k}, \ \ \ A_0=A_{x^k}y^k, \ \ \ A_{0j}=A_{x^ky^j}y^k.
\]
Suppose that $(A_{ij})$ is a positive definite tensor and $(A^{ij})$ denotes its inverse.   Then the following  hold
\begin{eqnarray}
g_{ij}= \frac{A^{\frac{2}{m}-2}}{m^2}[mAA_{ij}+(2-m)A_iA_j],\ \ \ g^{ij}= A^{-\frac{2}{m}}[mAA^{ij}+\frac{m-2}{m-1}y^iy^j], \ \ \\
y^iA_i=mA, \ \ y^iA_{ij}=(m-1)A_j,\ \ y_i=\frac{1}{m}A^{\frac{2}{m}-1}A_i,\ \  A^{ij}A_i=\frac{1}{m-1}y^j, \ \ A_iA_jA^{ij}=\frac{m}{m-1}A. \ \
\end{eqnarray}
Let $(M, F)$ be a Finsler manifold. The second derivatives of ${1\over 2} F_x^2$ at $y\in T_xM_{0}$ is an inner product ${\bf g}_y$ on $T_xM$.  The third order derivatives of ${1\over 2} F_x^2$ at  $y\in T_xM_0$ is a symmetric trilinear forms ${\bf C}_y$ on $T_xM$. We call ${\bf g}_y$ and ${\bf C}_y$ the  fundamental form and  the Cartan torsion, respectively.  The rate of change of the Cartan torsion along geodesics is the  Landsberg curvature  ${\bf L}_y$ on $T_xM$ for any $y\in T_xM_0$. $F$ is said to be   Landsbergian if ${\bf L}=0$. The quotient ${\bf L}/{\bf C}$ is regarded as the relative rate of change of Cartan torsion ${\bf C}$ along Finslerian geodesics.  Then $F$  is said to be isotropic Landsberg metric if ${\bf L}=cF \bf C$, where $c=c(x)$ is a scalar function on $M$. In this paper, we prove the following.

\bigskip

\begin{thm}\label{MainTHM1}
Let $(M,F)$ be an $n$-dimensional $m$-th root Finsler manifold. Suppose that $F$ is a non-Riemannian isotropic Landsberg metric. Then $F$ reduces to a Landsberg  metric.
\end{thm}

\bigskip

Let $F$ be a Finsler metric on a manifold $M$. The geodesics of $F$ are characterized locally by the equations $ {d^2 x^i\over dt^2} + 2 G^i( x, {dx\over dt})=0$, where $G^i = {1\over 4} g^{ik} \Big \{ 2 {\pa g_{pk}\over \pa x^q}
- {\pa g_{pq}\over \pa x^k} \Big \} y^py^q$ are coefficients of  the spray associated with $F$.  A Finsler metric $F$ is  called a Berwald metric if $G^i = {1\over 2} \Gamma^i_{jk}(x)y^jy^k$ are quadratic in $y\in T_xM$  for any $x\in M$. Taking a trace of Berwald curvature gives rise the mean Berwald curvature ${\bf E}$.  In \cite{AZ1}, Akbar-Zadeh introduces the non-Riemannian quantity $\bf{H}$ which is obtained from the mean Berwald curvature by the covariant horizontal differentiation along geodesics. He proves that for a Finsler manifold of scalar flag curvature ${\bf K}$ with dimension $n\geq3$, ${\bf K}=constant$\ \ if and only if ${\bf{H}}=0$. It is remarkable that, the Riemann curvature $R_y= R^i_{\ k}  dx^k \otimes \pxi|_x : T_xM \to T_xM$ is a family of linear maps on tangent spaces, defined by
\[
R^i_{\ k} = 2 {\pa G^i\over \pa x^k}-y^j{\pa^2 G^i\over \pa x^j\pa y^k}
+2G^j {\pa^2 G^i \over \pa y^j \pa y^k} - {\pa G^i \over \pa y^j}{\pa G^j \over \pa y^k}.
\]
A Finsler metric $F$ is said to be of scalar curvature if there is a scalar function ${\bf K}={\bf K}(x, y)$ such that $R^i_{\ k}= {\bf K}(x, y)F^2h^i_k$. If ${\bf K}=constant$, then $F$ is called of constant flag curvature.

A Finsler metric is called of isotropic ${\bf H}$-curvature if $H_{ij}= \frac{n+1}{2F}\theta h_{ij}$,  for some 1-form $\theta$ on $M$, where $h_{ij}$ is the angular metric. It is remarkable that in \cite{NST},  Z. Shen with the authors  prove that every Finsler metric of scalar flag curvature ${\bf K}$ and of isotropic  ${\bf H}$-curvature has almost isotropic flag curvature, i.e., the flag curvature is in the form ${\bf K}=\frac{3\theta}{F}+\sigma$,  for some scalar function $\sigma$ on $M$.

\bigskip

\begin{thm}\label{MainTHM3}Let $(M,F)$ be an $n$-dimensional $m$-th root manifold with $n\geq 2$. Suppose that $F$ has isotropic ${\bf H}$-curvature. Then ${\bf H}=0$.
\end{thm}
\section{Proof of the Main Theorems}
\begin{lem}\label{YY}{\rm (\cite{YY})}
\emph{Let  $F$ be an $m$-th root Finsler metric on an open subset $U\subset \mathbb{R}^n$. Then the spray coefficients of $F$ are given by}
\be
G^i=\frac{1}{2}(A_{0j}-A_{x^j})A^{ij}.\label{YY1}
\ee
\end{lem}

\noindent {\bf Proof of Theorem \ref{MainTHM1}}: Let $F=\sqrt[m]{A}$ be an $m$-th root isotropic Landsberg metric, i.e., $L_{ijk}=cFC_{ijk}$, where $c=c(x)$ is a scalar function on $M$. The Cartan tensor of $F$ is given by the following
\be \label{Cartan}
C_{ijk}=\frac{1}{m} A^{\frac{2}{m}-3}\Big[A^2A_{ijk}+(\frac{2}{m}-1)(\frac{2}{m}-2)A_iA_jA_k+(\frac{2}{m}-1)A\{A_iA_{jk}+A_jA_{ki}+A_kA_{ij}\}   \Big].
\ee
Since $L_{ijk}=-\frac{1}{2}y_sG^s_{\ y^iy^jy^k}$, then we have $L_{ijk}=-\frac{1}{2m}A^{\frac{2}{m}-1}A_sG^s_{\ y^iy^jy^k}$. Therefore, we have
\be
A_sG^s_{\ y^iy^jy^k}=-2c A^{\frac{1}{m}-2}\Big[A^2A_{ijk}+(\frac{2}{m}-1)\Big\{(\frac{2}{m}-2)A_iA_jA_k+A\{A_iA_{jk}+A_jA_{ki}+A_kA_{ij}\} \Big\}  \Big]. \label{L1}
\ee
By Lemma \ref{YY}, the left hand side of (\ref{L1}) is a rational function in $y$, while its right hand side is an irrational  function in $y$.  Thus, $c=0$ or $A$ satisfies the following PDEs
\be \label{eq 0}
A^2A_{ijk}+(\frac{2}{m}-1)(\frac{2}{m}-2)A_iA_jA_k+(\frac{2}{m}-1)A\{A_iA_{jk}+A_jA_{ki}+A_kA_{ij}\}=0
\ee
Plugging (\ref{eq 0}) into (\ref{Cartan}) implies that $C_{ijk}=0$. Hence, by Deicke's theorem, $F$ is Riemannain metric, which
contradicts our assumption. Therefore, $c=0$. This completes the proof.
\qed

\bigskip

\noindent {\bf Proof of Theorem \ref{MainTHM3}}: Let $F=\sqrt[m]{A}$  be of isotropic ${\bf H}$-curvature, i.e.,
\be \label{IMB}
{H}_{ij}= \frac{n+1}{2F} \theta h_{ij},
\ee where $\theta$ is a 1-form on $M$. The angular metric $h_{ij}=g_{ij}-F^{2}y_iy_j$ is given by the following
\be\label{angular}
h_{ij}= \frac{1}{m^2}[mAA_{ij}+(1-m)A_iA_j]A^{\frac{2}{m}-2}.
\ee
Plugging (\ref{angular}) into (\ref{IMB}), we get
\be \label{E}
H_{ij}=\frac{(n+1)\theta}{2m^2}\ [mAA_{ij}+(1-m)A_iA_j]A^{\frac{1}{m}-2}.
\ee
By (\ref{YY1}), one can see that $H_{ij}$ is rational with respect to $y$. Thus, (\ref{E}) implies that
$\theta=0$ or
\be\label{E1}
mAA_{ij}+(1-m)A_iA_j=0.
\ee
By (\ref{angular}) and (\ref{E1}), we conclude that $h_{ij}=0$, which is impossible. Hence $\theta=0$ and  then $H_{ij}=0$.
\qed

\bigskip

By Schur Lemma, Theorem \ref{MainTHM3} and Theorem 1.1 of \cite{NST}, we have the following.\\

\begin{cor}
Let $(M,F^n)$ be an n-dimensional m-th root Finsler manifold of scalar flag curvature ${\bf K}$ with $n\geq 3$.  Suppose that the flag curvature is given by ${\bf K}=\frac{3\theta}{F}+\sigma$, where $\theta$ is a 1-form and
$\sigma=\sigma(x)$ is a scalar function on M. Then ${\bf K}=0$.
\end{cor}

\bigskip
\noindent
Akbar Tayebi\\
Faculty  of Science, Department of Mathematics\\
Qom University\\
Qom. Iran\\
Email: akbar.tayebi@gmail.com
\bigskip
\noindent
Behzad Najafi\\
Faculty  of Science, Department of Mathematics\\
Shahed University\\
Tehran. Iran\\
najafi@shahed.ac.ir
\bigskip

\noindent

\end{document}